\documentclass[11pt,fleqn,twoside]{article}
\usepackage{amsfonts,amssymb,latexsym}
\makeatletter
\newcommand{\prava}[1]{\small\it
\begin{flushleft}
Copyright \copyright \ 1999 by  #1
\end{flushleft}}

\newcommand{\name}[1]{\begin{flushleft}
                       \LARGE \bf #1
                       \end{flushleft}\vspace{-3mm}}

\newcommand{\Author}[1]{\begin{flushleft}
                       \it #1 \end{flushleft}}

\newcommand{\Adress}[1]{\begin{flushleft}
                       \it #1 \end{flushleft}}

\newcommand{\Date}[1]{\begin{flushleft}
                      \small  \it #1 \end{flushleft}}

\newcommand{\ehkol}{Author \ name}
\newcommand{\ohkol}{Article \ name}
\renewcommand{\@evenhead}{
\hspace*{-3pt}\raisebox{-15pt}[\headheight][0pt]{\vbox{\hbox to \textwidth 
{\thepage \hfil \ehkol}\vskip4pt \hrule}}}
\renewcommand{\@oddhead}{
\hspace*{-3pt}\raisebox{-15pt}[\headheight][0pt]{\vbox{\hbox to \textwidth 
{\ohkol \hfil \thepage}\vskip4pt\hrule}}}
\renewcommand{\@evenfoot}{}
\renewcommand{\@oddfoot}{}

\newcommand{\be}{\begin{equation}}
\newcommand{\ee}{\end{equation}}
\newcommand{\ba}{\hspace*{-5pt}\begin{array}}
\newcommand{\ea}{\end{array}}
\newcommand{\p}{\partial}
\newcommand{\ds}{\displaystyle}
\makeatother

\def\dd{{\rm d}}

\begin{document}
\thispagestyle{empty}
\setcounter{page}{355}
\renewcommand{\ehkol}{C. Sophocleous}
\renewcommand{\ohkol}{Continuous and Discrete Transformations of 
a Porous Medium Equation}

\begin{flushleft}
\footnotesize \sf Journal of Nonlinear Mathematical Physics \qquad
1999, V.6, N~4, \pageref{sophocleous-fp}--\pageref{sophocleous-lp}.
\hfill {\sc Letter}
\end{flushleft}

\vspace{-5mm}

\renewcommand{\footnoterule}{}
{\renewcommand{\thefootnote}{}
 \footnote{\prava{C. Sophocleous}}}

\name{Continuous and Discrete Transformations of\\
a One-Dimensional Porous Medium Equation}\label{sophocleous-fp}

\Author{Christodoulos SOPHOCLEOUS}

\Adress{Department of Math. and Stat.,
University of Cyprus, CY-1678 Nicosia, Cyprus\\
E-mail: christod@olympia.mas.ucy.ac.cy}

\Date{Received March 07, 1999; Revised July  09, 1999;
Accepted July 12, 1999}

\begin{abstract}
\noindent
We consider the one-dimensional porous medium equation 
$u_t=\left (u^nu_x \right )_x+\frac{\mu}{x}u^nu_x$. We derive 
point transformations of a general class
that map this equation into itself or into equations
of a similar class. In some cases this porous medium equation is connected
with well known equations. With the introduction of a new dependent variable 
this partial dif\/ferential equation
can be equivalently written as a system of two equations. Point 
transformations are also sought for this auxiliary system. It turns out that
in addition to the continuous point transformations that may be derived
by Lie's method, a number of discrete transformations are obtained. In some
cases the point transformations which are presented here for the
single equation and for the auxiliary system form cyclic groups of
f\/inite order. 
\end{abstract}

\section{Introduction}

\renewcommand{\theequation}{\thesection.\arabic{equation}}
\setcounter{equation}{0}

Probably the most useful point transformations of partial dif\/ferential 
equations (pdes) are those which form a continuous (Lie) group of 
transformations, each member of which leaves an equation invariant.
In addition to possessing continuous groups of symmetries, many pdes also
possess discrete symmetries. However, the classical method of deriving the
continuous symmetries, which is based on the linear nature of inf\/initesimal
transformations, is not directly applicable to the derivation of groups of
discrete symmetries. Such example is given by Kingston and Sophocleous~[1]
who found that the reciprocal point transformation (double application
gives the identity transformation) $x' =x/t$, $t' =1/t$, $u' =-ut+x$ leaves
the Burger-type equation $u_t+uu_x+(f(t)-f(1/t))u_{xx}=0$ invariant, which is
a symmetry additional to the Lie point symmetries obtained from the classical 
approach.

Point transformations may also be used to link a pde
with an equation of a dif\/ferent form. For example, it is sometimes
possible to connect an equation to a canonical form for which there
is an established theory. Such cases are given in the present paper.

For these reasons it is advantageous to study point transformations directly
in f\/inite forms with the ultimate dual aims of f\/inding the complete set of 
point transformation symmetries of pdes and also discovering new links between
dif\/ferent equations. A good example, of the approach which is used here,
is the work by Sophocleous and Kingston who classify the point transformations
of four common classes of one-dimensional non-linear wave equations~[2].

We consider the one-dimensional porous medium equation of the form
\begin{equation}
u_t=\left (u^nu_x \right )_x+\frac{\mu}{x}u^nu_x .
\end{equation}
We introduce point transformations of a general class and investigate 
when such transformations map Eq.(1.1) into itself. In addition we classify
point transformations that link Eq.(1.1) to dif\/ferent equations but of the
same class, where in some cases it turns out to be well known equations.
A number of the point transformations obtained form
cyclic groups of f\/inite order. 

The next section is devoted to the above analysis, where f\/irstly it makes use 
of a result obtained by Kingston and Sophocleous~[3] which is 
summarised in the following theorem:

\medskip
\noindent
{\bf Theorem.} {\it The point transformation $x' =P(x,t,u)$,
$t' =Q(x,t,u)$, $u' =R(x,t,u)$ transforms the pde
\[
u'_{t'}=H'(x',t',u',u'_1,\ldots,u'_n)
\]
to the pde
\[
u_t=H(x,t,u,u_1,\ldots,u_n)
\]
where $u_n=\frac{\p ^nu}{\p x^n}$, $u'_n=\frac{\p^n u'}{\p x'^n}$,
$n\ge 2$ and $H$ and $H'$ are polynomials (non-negative integral powers) in
$u_1,\ldots,u_n$ and $u'_1,\ldots,u'_n$ respectively, if and only if
\[
P=P(x,t),\qquad Q=Q(t)
\]
and the functions $H$ and $H'$ satisfy the relation
\begin{equation}
P_xR_uH=P_xQ_tH'+P_tR_x+P_tR_uu_x-P_xR_t .
\end{equation}}

The proof of this theorem can be found in ref.~[3].

Equation (1.1) can be written in a conserved form and therefore with an 
introduction of a new dependent variable may be represented by a system
of two pdes. In the Sections~3 and~4 we search for point transformations of 
a general class that map this system into itself or to a system of the same 
class.
In some cases, such point transformations of this auxiliary system lead 
to non-local transformations of Eq.(1.1).

Equation (1.1) is a special case of the generalised porous medium equation
\begin{equation}
u_t=\left (u^nu_x \right )_x+f(x)u^su_x+g(x)u^m .
\end{equation}
The classical point symmetries have been presented by Gandarias [4].
Further study, along the lines of this paper, of Eq.(1.3) may therefore
be useful.

\section{Point transformations}

\setcounter{equation}{0}

We consider point transformations of the form
\begin{equation}
x' =P(x,t,u), \qquad t' =Q(x,t,u), \qquad u' =R(x,t,u)
\end{equation}
which transform the pde
\begin{equation}
u'_{t'}=\left ( f(u' )u'_{x'} \right )_{x'}+g(x' ,u' )u'_{x'}
\end{equation}
to the pde
\begin{equation}
u_t=\left (F(u)u_x \right )_x+ G(x,u) u_x,
\end{equation}
where $F$ not a constant. We assume that point transformation (2.1) is
non-degenerate in the sense that
\[
\frac{\p (P,Q,R)}{\p (x,t,u)}\ne 0 \qquad \mbox{and} \qquad 
\frac{\p (P(x,t,u(x,t),Q(x,t,u(x,t)))}{\p (x,t)} \ne 0.
\]

Firstly, we note that each linear transformation 
of the form
\begin{equation}
x' =c_1x+c_2, \qquad t' =c_3t+c_4, \qquad u' =c_5u+c_6,
\end{equation} 
where $c_i$ are constants (possibly complex), $c_1c_3c_5 \ne 0$,
achieves this purpose. That is, Eq.(2.2) is transformed to Eq.(2.3)
with the function $F(u)$ def\/ined in terms of $f(u' )$ and $G(x,u)$ in
terms of $g(x',u')$ by
\begin{equation}
F(u)=c_3c_1^{-2}f(c_5u+c_6), \qquad G(x,u)=c_3c_1^{-1}g(c_1x+c_2,c_5u+c_6),
\end{equation}
respectively. Hence, point transformations of the form (2.4) may be f\/iltered 
out of equations arising in the discussion, without loss of generality,
with the understanding that all point transformations obtained may be 
augmented by Eqs.(2.4) and (2.5).

We now turn into the special case of Eq.(2.3), where $F(u)=u^n$ and
$G(x,u)=\frac{\mu}{x}u^n$,
\begin{equation}
u_t=\left (u^nu_x \right )_x+\frac{\mu}{x}u^nu_x.
\end{equation}
The Lie point symmetries of the equation 
$u_t=\left (u^n\right )_{xx}+\frac \mu x 
\left (u^n \right )_x$, which is equivalent to Eq.(2.6) with the exception
of the case where $n=-1$, have been classif\/ied by Gandarias et~al~[5]. 
The symmetries of Eq.(2.6) are summarised in the following table, where
$(X,T,U)$ are the components of the symmetry generator 
$X\frac{\p~}{\p x}+T\frac{\p~}{\p t}+U\frac{\p~}{\p u}$:

\medskip

\begin{equation}
\begin{tabular}{lccl}
\hline
&$n$&$\mu$&$(X,T,U)$ \\
\hline
&&&\\[-4mm]
$\Gamma_1$&$\ne 0$&arbitrary&$(0,1,0)$ \\
$\Gamma_2$&&&$(x,2t,0)$ \\
$\Gamma_3$&&&$(0,nt,- u)$ \\[1mm]
\hline
&&&\\[-3mm]
$\Gamma_4$&$\ne -2,-1,0$&$=\frac{3n+4}{n+2}$&$\left ((n+2)x^{-\frac{n}{n+2}},
0,-2x^{-\frac{2(n+1)}{n+2}}u\right )$ \\[3mm]
\hline
&&&\\[-3mm]
$\Gamma_5$&$=-1$&$=1$&$(x\ln x,0,2(\ln x-1)u)$ \\[1mm]
\hline
\end{tabular}
\end{equation}

\medskip

We search for point transformations of the form (2.1) that map
Eq.(2.2) into Eq.(2.6). We note that both of these equations are such
that the Theorem, given in the Introduction, can be applied.
Therefore we have $Q=Q(t)$ and $P=P(x,t)$. Also we substitute
$H=u^nu_{xx}+nu^{n-1}u_x^2+\frac \mu x u^nu_x$,
$H'=f(u')u'_{x'x'}+\frac{\dd f}{\dd u'}{u'_{x'}}^2+g(x',u')u'_{x'}$ 
into Eq.(1.2) and with the use of the expressions of $u'_{x'}$,
$u'_{x'x'}$ [3] the relation (1.2) becomes
\begin{equation}
\ba{l}
\ds xuP_xR_u\left (u^n P_x^2-Q_tf \right )u_{xx}
\vspace{3mm}\\
\ds  + xP_x\left (nu^nP_x^2R_u-uQ_tR_{uu}f-uQ_tR_u^2
\frac{\dd f}{\dd u'} \right )u_x^2 + 
 u\Bigl (xQ_tP_{xx}R_uf-xP_tP_x^2R_u
\vspace{3mm}\\
\ds  +\mu u^nP_x^3R_u-xQ_tP_x^2R_ug-
2xQ_tP_xR_{xu}f-2xQ_tP_xR_xR_u \frac{\dd f}{\dd u'} \Bigr )u_x 
\vspace{3mm}\\
\ds  +xu\left
(Q_tP_{xx}R_xf-P_tP_x^2R_x+P_x^3R_t-Q_tP_x^2R_xg-Q_tP_xR_{xx}f- 
Q_tP_xR_x^2 \frac{\dd f}{\dd u'} \right )=0.
\ea\hspace{-19.17pt}
\end{equation}
We set, successively, the coef\/f\/icients of $u_{xx}$, $u_x^2$,
$u_x$ and the term 
independent of $u_x$ and $u_{xx}$ in Eq.(2.8) equal to zero to obtain four
identities. These identities enable the desired point transformations to be
derived and ultimately impose restrictions on the functional forms of
$Q(t)$, $P(x,t)$, $R(x,t,u)$, $f(u' )$ and $g(x',u')$. For example, from
the f\/irst two we deduce that $f=\frac{P_x^2}{Q_t}u^n$ and also that $R$ is 
linear in $u$. From these we conclude that $f(u' )$ is linear in ${u'}^n$, 
but using the transformation (2.4) we can simply write $f(u' )={u'}^n$.

In the following analysis we omit any further calculations. It turns out that
this analysis can be split into two cases:(1) $n=-1$; (2) $n \ne -1$.
Firstly, we point out that the discrete transformations which
correspond to the symmetries $\Gamma_1$--$\Gamma_3$ admitted by Eq.(2.6)
are of the form (2.4).

{\bf Case 1}: ($n=-1$) Equation (2.6) takes the form
\begin{equation}
u_t=\left (\frac{u_x}{u} \right )_x+\frac{\mu}{xu}u_x.
\end{equation}
If $\mu =1$, Eq.(2.9) is mapped into itself by the transformation
\begin{equation}
x' =x^{k+1}, \qquad t' =t, \qquad u' =\frac {1}{(k+1)^2}x^{-2k}u .
\end{equation}
This latter transformation is the continuous symmetry which is
represented by the symmetry $\Gamma_5=x\ln x\frac{\dd~}{\dd x}+2(\ln
x-1)u\frac{\dd~}{\dd u}$ admitted by Eq.(2.9).

If we apply transformation (2.10) $N$-times we obtain
\begin{equation}
x^{(N)}=x^{(k+1)^N}, \qquad
t^{(N)}=t, \qquad u^{(N)}=\frac{1}{(k+1)^{2N}}x^{-2(k+1)^N+2}u.
\end{equation}
We note that this point transformation forms a cyclic group of order $N$ if
$k$ is a root of the equation
\begin{equation}
(k+1)^N-1=0.
\end{equation}
For example, if $N=2$ then $k=-2$ (we omit the obvious solution $k=0$) and
we get the transformation $x' =1/x$, $t' =t$, $u' =x^4u$ which forms a cyclic 
group of order 2. In fact, it is trivial to note that $k=0$ and $k=-2$ are 
the only real roots of equation (2.12). If we introduce complex roots 
of Eq.(2.12), then we may construct cyclic groups of any f\/inite order. 
For example, the point transformation $x' =x^{\rm i}$, $t' =t$, 
$u' =-x^{2-2{\rm i}}u$ forms a cyclic group of order 4.

Furthermore identity (2.8) produces the following three equations that mapped
into Eq.(2.9) ($\mu =1$) by the corresponding discrete transformations:
\be
\ds u'_{t'}= \left (\frac{u'_{x'}}{u'} \right )_{x'};
\qquad  
x' =k\ln |x|, \quad t' =t, \quad u' =\frac{1}{k^2}ux^2;
\ee
\be
\ds u'_{t'}= \left (\frac{u'_{x'}}{u'} \right )_{x'}+
\left (\frac{1}{x'u'} +\frac 12 x' \right )u'_{x'}; \qquad\!
\ds x' =\frac{x^{k}}{\sqrt{t}}, \quad t' =\ln t, \quad 
u' =\frac{1}{k^2}ux^{2(1-k)};
\ee
\be
\ds u'_{t'}= \left (\frac{u'_{x'}}{u'} \right )_{x'}+
\frac 12 x' u'_{x'}; \qquad
\ds x' =\frac{k\ln |x|}{\sqrt{t}}, \quad t' =\ln t,
\quad u' =\frac{1}{k^2}ux^2.
\ee

We note from result (2.13) that Eq.(2.9) is connected with a well known
non\-li\-near dif\-fusion equation. Therefore any solution of this dif\/fusion
equation can be transformed into a solution of Eq.(2.9) using (2.13).
For example, point transformation (2.13) maps the solution 
$u' = 2t' \left ({x'}^2+4{t'}^2 \right )^{-1}$ of
$u'_{t'}=(u'_{x'}/u')_{x'}$ into the solution 
$u=2t\left[x^2(\ln^2|x|+4t^2)\right ]^{-1}$ of Eq. (2.9). 

Finally, if $\mu =0$, the following equation with the corresponding
point transformation
\begin{equation}
u'_{t'}= \left (\frac{u'_{x'}}{u'} \right )_{x'}+
\left (\frac{1}{x'u'} +\frac 12 x' \right )u'_{x'}; \qquad
x' =\frac{{\rm e}^x}{\sqrt{t}}, \quad t' =\ln t , \quad u' =u{\rm e}^{-2x}
\end{equation}
is mapped into Eq.(2.9) which has the form of the nonlinear
dif\/fusion equation $u_t=(u_x/u)_x$.
Result (2.16) can also be obtained by combining the results
(2.13) and (2.14).

{\bf Case 2}: ($n\ne -1$) The continuous symmetry which corresponds to the 
inf\/initesimal symmetry $\Gamma_4=(n+2)x^{-\frac{n}{n+2}}\frac{\p~}{\p x}-2u
x^{-\frac{2(n+1)}{(n+2)}}\frac{\p~}{\p u}$ admitted by Eq.(2.6) is 
the following
\begin{equation}
x' =\left ( x^{\frac{2n+2}{n+2}} + C \right )^{\frac{n+2}{2n+2}},
\qquad t'=t, \qquad 
u' =ux^{\frac{2}{n+2}}\left ( x^{\frac{2n+2}{n+2}}+C \right
)^{-\frac{1}{n+1}}, 
\end{equation}
where also $n\ne -2,0$ and $\mu =\frac {3n+4}{n+2}$.

In addition to the above result, identity (2.8) produces three more discrete
transformations. When $\mu\ne -\frac{n+2}{n}$, the point transformation
\begin{equation}
x' =\frac{2n+2}{\mu n+n+2}x^{\frac{\mu n+n+2}{2n+2}},
\qquad t' =t, \qquad u' =x^{\frac{\mu -1}{n+1}}u
\end{equation}
maps the equation
\begin{equation}
u'_{t'}=\left ({u'}^nu'_{x'} \right )_{x'}+\lambda 
\frac{{u'}^n}{x'}u'_{x'}
\end{equation}
into Eq.(2.6), where $\lambda =\frac{3n+4-n\mu -2\mu}
{\mu n+n+2}$. If $\lambda =\mu =-\frac{3n+4}{n}$ then the discrete point
transformation $x' =\frac{1}{x}$, $t' =t$, $u' =x^{-\frac 4n}u$,
which forms a cyclic group of order 2, maps Eq.(2.6) into itself. We
observe that $\mu =\frac{3n+4}{n+2}$ 
implies $\lambda =0$ and Eq.(2.19) becomes the standard nonlinear dif\/fusion 
equation which is transformed to Eq.(2.6) ($n \ne -2,-1$) by the
transformation 
\[
x' =\frac{n+2}{2n+2}x^{\frac{2n+2}{n+2}}, 
\qquad t'=t, \qquad u' =x^{\frac {2}{n+2}}u .
\]
Also $\mu =0$ implies $\lambda =\frac{3n+4}{n+2}$ and we obtain the inverse
of this latter point transformation which maps Eq.(2.19) into the
nonlinear dif\/fusion equation.

When $\mu =-\frac{n+2}{n}$, the point transformation
\begin{equation}
x' =\ln |x|,\qquad t' =t, \qquad u' =x^{-\frac 2n}u
\end{equation}
maps the equation
\begin{equation}
u'_{t'}=\left ({u'}^nu'_{x'} \right )_{x'}+2\frac{n+1}{n}
{u'}^nu'_{x'}
\end{equation}
into Eq.(2.6). Eq.(2.21) is the Boussinesq equation of hydrology involved
in various f\/ields of petroleum technology and ground water hydrology.
We observe that when $n=-2$ ($\mu =0$) Eq.(2.21) is mapped by (2.20)
into the well known nonlinear dif\/fusion equation $u_t=\left (u^{-2}u_x
\right )_x$. This dif\/fusion equation has a number of properties. For
example, admits Lie-B\"acklund transformations and there exists a 
transformation that maps it into the linear heat equation $u_t=u_{xx}$~[6].
Also we point out that if $n=-1$ we obtain the result (2.13) of the case 1.

Finally, when $\mu \ne -\frac{n+2}{n}$, the point transformation
\begin{equation}
x' =\frac{2n+2}{\mu n+n+2}\frac{x^{\frac{\mu n+n+2}{2n+2}}}{\sqrt{t}}, 
\qquad t' =\ln t, \qquad u' =x^{\frac{\mu -1}{n+1}}u
\end{equation}
maps the equation
\begin{equation} 
u'_{t'}=\left ({u'}^nu'_{x'} \right )_{x'}+\left (\lambda 
\frac{{u'}^n}{x'}+\frac 12 x' \right )u'_{x'}
\end{equation}
into Eq.(2.6), where  $\lambda =\frac{3n+4-n\mu -2\mu}{\mu n+n+2}$. 
In the case where $\mu =-\frac{n+2}{n}$, identity (2.8) gives
$n=-1$ ($\mu =1$, $\lambda =0$) and we simply reproduce the result (2.15).

\section{Potential transformations}

\setcounter{equation}{0}

If we introduce the potential $v$, Eq.(2.6) can be written as a
system of two pdes,
\begin{equation}
v_x=ux, \qquad v_t=xu^nu_x+\frac{\mu-1}{n+1}u^{n+1},
\end{equation}
if $n\ne -1$ and as
\begin{equation}
v_x=ux, \qquad v_t=\frac xu u_x+(\mu -1)\ln |u|,
\end{equation}
if $n=-1$.

Similarly, if we write 
\begin{equation}
g(x',u')=\frac{f(u' )+h'(u' )}{x'}, 
\end{equation}
where $h'(u' )=\frac{\dd h}{\dd u'}$, then Eq.(2.2) can also be written 
as a system of two pdes
\begin{equation}
v'_{x'} =x'u', \qquad v'_{t'}=x' f(u')u'_{x'}+h(u').
\end{equation}
We also consider the point transformation
\begin{equation}
x' =P(x,t,u,v), \quad t' =Q(x,t,u,v), \quad
u' =R(x,t,u,v), \quad v' =S(x,t,u,v)
\end{equation}
relating $x$, $t$, $u(x,t)$, $v(x,t)$ and $x'$, $t'$, $u'(x',t')$,
$v'(x',t')$ and we assume that it is non-degenerate.

It is well known that in many cases pdes which can be written in a conserved 
form admit nonlocal symmetries, known as {\it potential symmetries}~[7]. For
example, in order to f\/ind potential symmetries for Eq.(2.6) we search for 
point symmetries for the system (3.1) (and (3.2)). If at least one of the
inf\/initesimal generators of the variable $x$, $t$, $u$ depends on the 
potential variable $v$, then the point symmetry of the auxiliary system (3.1)
is also a potential symmetry of the Eq.(2.6). Otherwise, the point symmetry
of (3.1) projects onto a point symmetry of Eq.(2.6).

The symmetry analysis of the system (3.1) was carried out in the ref.~[5]. 
It was shown that only one symmetry of (3.1) induces a potential
symmetry admitted by Eq.(2.6). The rest of the symmetries project
onto point symmetries of Eq.(2.6). The potential symmetry occurs if
$n=-2$ and $\mu =-\frac 12$ and is given by
\begin{equation}
2xv\frac{\p~}{\p x}-2(x^2u^2+uv)\frac{\p~}{\p u}+v^2\frac{\p~}{\p v} .
\end{equation}
Symmetry (3.6) generates the one-parameter continuous group of
transformations \begin{equation}
x' =\frac{x}{(1-\epsilon v)^2}, \qquad t' =t, \qquad
u' =\frac{u(1-\epsilon v)^3}{2\epsilon x^2u+1-\epsilon v},
\qquad v' =\frac{v}{1-\epsilon v}.
\end{equation}

Our goal in this section is to derive point transformations of the class
(3.5) which map Eqs.(3.4) into Eqs.(3.1) (and Eqs.(3.4) into Eqs.(3.2)). 
Using (3.5) and the
forms of $u'_{x'}$, $u'_{t'}$, $v'_{x'}$, $v'_{t'}$ in terms of
$u_x$, $u_t$, $v_x$, $v_t$ and derivatives of $P$, $Q$, $R$, $S$ [8]
and eliminating $v_x$ and $v_t$ from Eqs.(3.1) (or Eqs.(3.2)),
Eqs.(3.4) become two identities of the form
\begin{equation}
E_1(x,t,u,v,u_x,u_t)=0, \qquad E_2(x,t,u,v,u_x,u_t)=0,
\end{equation}
where $x$, $t$, $u$, $v$, $u_x$ and $u_t$ are regarded as independent
variables and $E_1$, $E_2$ are explicit polynomials in $u_x$ and $u_t$.

Now, the coef\/f\/icients of $u_x^2$, $u_x$, $u_t$ in $E_1=0$ and
$u_x^2$, $u_t$ in $E_2=0$, which must be identically equal to zero, give
$P_u=Q_x=Q_u=Q_v=S_u=0$. Therefore, the point transformation (3.5) can be
written in the simplif\/ied form
\begin{equation}
x' =P(x,t,v), \qquad t' =Q(t), \qquad u' =R(x,t,u,v), \qquad v' =S(x,t,v).
\end{equation}
To ensure that the point transformations are non-degenerate we require to have
\begin{equation}
Q_tR_u(P_xS_v-P_vS_x) \ne 0.
\end{equation}
In addition, we need to have
\begin{equation}
P_v^2+R_v^2 \ne 0
\end{equation}
because otherwise the point transformations derived are equivalent to point 
transformations which connect Eq.(2.2) and Eq.(2.6). Transformations that 
satisfy condition (3.11) shall be called {\it potential transformations}.

The restricted forms (3.9) simplify the identities (3.8). For the convenience
of the reader we omit any further calculations. We present the three
potential transformations obtained: (1) $n=-1$; (2) $n=-2$; (3) $n\ne
-2,-1,0$.

{\bf Case 1}: If $n=-1$ then the point transformation
\begin{equation}
x' ={\rm e}^{t+\frac v2}, \qquad t' =t, \qquad 
u' =\frac{4\epsilon {\rm e}^{-(v+2t)}}{ux^2}, \qquad 
v' =2(\epsilon\ln x-t),
\end{equation}
transforms the pdes
\begin{equation}
v'_{x'}=x'u', \qquad v'_{t'}=\epsilon \frac{x'}{u'}u'_{x'}+2(\epsilon -1),
\end{equation}
where $\epsilon$ is a constant, to the pdes (3.2) where $\mu =1$. 
Application of the point transformation (3.12) 2N-times gives
\begin{equation}
{x'}^{(2N)}=x^{\epsilon^N},
\quad {t'}^{(2N)}=t, \quad {u'}^{(2N)}=ux^{2(1-\epsilon^N)}, \quad
{v'}^{(2N)}=\epsilon^Nv+2(\epsilon^N-1)t .
\end{equation}
Therefore transformation (3.12) forms a cyclic group of order $2N$ if
$\epsilon$ is a root of the equation $\epsilon^N=1$. Clearly, the only two
real roots are $\epsilon =1$ and $\epsilon =-1$ where in these cases
transformation (3.12) forms a cyclic group of order 2 and order 4,
respectively. If we allow $\epsilon$ to be complex, then we may construct
cyclic groups of any even order. For example, if $\epsilon ={\rm i}$ 
then the point transformation $x' ={\rm e}^{t+v/2}$, $t' =t$,
$u' =\frac{4{\rm i}}{ux^2\exp(v+2t)}$, $v' =2({\rm i}\ln x-t)$ forms a cyclic
group of order 8.

If $\epsilon =1$ the reciprocal transformation (3.12) maps Eqs.(3.2) ($\mu =1$)
into itself and therefore we also observe that the point transformation
\[
x' ={\rm e}^{t+v/2}, \qquad t' =t, \qquad v' =2(\ln x-t)
\]
leaves the integrated form of Eq.(2.9) ($u=\frac{v_x}{x}$),
\[
v_t=\frac{x}{v_x}v_{xx}-1
\]
invariant.

{\bf Case 2}: If $n=-2$ and $\mu =-\frac 12$ then the point transformation
\begin{equation}
x' =-\frac{x}{v^2}, \qquad 
t' =t, \qquad u' =\frac{uv^3}{2x^2u-v}, \qquad v' =\frac{1}{v}
\end{equation}
maps the system (3.1) into itself. This latter
point transformation forms a cyclic group of order 2. From this result we 
deduce that the integrated form of Eq.(2.6),
\[
v_t=\left (\frac {x}{v_x} \right )^2v_{xx}+\frac 12 \frac {x}{v_x}
\]
remains invariant under the point transformation
$x' =-x/v^2$, $t' =t$, $v' =1/v$. 

\newpage

{\bf Case 3}: If $n \ne -2,-1,0$ and $\mu =\frac{3n+4}{n+2}$ then the point transformation
\begin{equation}
x' =v^{\frac{n}{2n+2}},~~t' =t, \qquad
u' =x^{-\frac{2}{n+2}}v^{\frac{1}{n+1}}u^{-1}, \qquad 
v' =\frac{n(n+2)}{4(n+1)^2} x^{\frac{2n+2}{n+2}},
\end{equation}
transforms the pdes
\begin{equation}
v'_{x'}=x'u', \qquad v'_{t'}=\lambda_1x'{u'}^{-(n+2)}u'_{x'}
+\lambda_2{u'}^{-(n+1)}
\end{equation}
to the pdes (3.1),
where $\lambda_1 =\frac{n^2}{4(n+1)^2}$ and
$\lambda_2=-\frac{n}{2(n+1)^2}$. If we use the integrated forms 
of Eq.(2.6) and Eq.(2.2) we deduce that the point transformation
\[
x' =v^{\frac{n}{2n+2}}, \qquad t' =t, \qquad 
v' =\frac{n(n+2)}{4(n+1)^2} x^{\frac{2n+2}{n+2}}
\]
maps the pde
\[
v'_{t'}=\lambda_1 \left (\frac{v'_{x'}}{x'} \right )^{-(n+2)}v'_{x'x'}
+(\lambda_2-\lambda_1) \left (\frac{v'_{x'}}{x'} \right )^{-(n+1)}
\]
into the pde
\[
v_t=\left (\frac{v_x}{x} \right )^{n}v_{xx}-\frac{n}{n+2}
\left (\frac{v_x}{x} \right )^{n+1}.
\]

\section{More potential transformations}

\setcounter{equation}{0}

We consider the standard nonlinear dif\/fusion equation
\begin{equation}
u'_{t'}=\left (f(u' )u'_{x'} \right )_{x'}
\end{equation}
which can be written, with the introduction of the potential variable $v$,
as a system of two pdes
\begin{equation}
v'_{x'}=u', \qquad v'_{t'}=f(u' )u'_{x'} .
\end{equation}
Here we present point transformations of the form (3.9) which map
Eqs.(4.2) into Eqs.(3.1) (or into Eqs.(3.2)).

{\bf Case 1}: If $n=-1$ the point transformation
\begin{equation}
x' =v+2t, \qquad t' =t, \qquad u' =\frac{c}{ux^2}+1,
\qquad v' =v+2t+c\ln x,
\end{equation}
transforms the system (4.2) with $f=\frac{c}{u'-1}$, $c$ a nonzero constant, 
into the system (3.2) where $\mu =1$.

{\bf Case 2}: If $n \ne -2,-1$ the point transformation
\begin{equation}
x' =v, \qquad t' =t, \qquad
u' =x^{-\frac{2}{n+2}}u^{-1}, \qquad 
v' = \frac{n+2}{2(n+1)}x^{\frac{2(n+1)}{n+2}},
\end{equation}
maps the system (4.2) with $f={u'}^{-(n+2)}$ into the system (3.1)
where $\mu =\frac{3n+4}{n+2}$. If we use the integrated forms of Eq.(2.6) 
and Eq.(4.1) we deduce that the point transformation $x' =v$, $t' =t$,
$v' =\frac{n+2}{2n+2}x^{\frac{2n+2}{n+2}}$ maps the pde
\[
v'_{t'}={v'_{x'}}^{-(n+2)}v'_{x'x'} 
\]
into the pde
\[
v_t=\left (\frac{v_x}{x} \right )^nv_{xx}-\frac{n}{n+2}
\left (\frac{v_x}{x} \right )^{n+1}.
\] 
We note that if $n=0$ we obtain the hodograph transformation $x' =v$,
$t' =t$,
$v' =x$ which connects the pde $v'_{t'}={v'_{x'}}^{-2}v'_{x'x'}$
and the linear heat equation $v_t=v_{xx}$. This latter transformation is
also a special case of the known general result where this later hodograph 
transformation connects the equations $v_t=F(v_x)v_{xx}$ and
$v_t=(v_x)^{-2}F(1/v_x)v_{xx}$ [9].

Furthermore point transformation (4.4) leads to the contact transformation
\begin{equation}
\dd x' =ux\dd x + (xu)_x \dd t , \qquad
\dd t' =\dd t, \qquad u' =x^{-\frac{2}{n+2}} u^{-1}
\end{equation}
which connects the nonlinear dif\/fusion equation 
$u'_{t'}=\left ({u'}^{-(n+2)}u'_{x'} \right )_{x'}$ and Eq.(2.6), 
where $\mu =\frac{3n+4}{n+2}$. 
We observe that if $n=0$, transformation (4.5)
maps the well known nonlinear dif\/fusion equation $u'_{t'}=
\left ({u'}^{-2}u'_{x'} \right )_{x'}$ into the linear pde
$u_t=u_{xx}+\frac{2}{x}u_x$.

\subsection*{Acknowledgements}

The author would like to thank Prof. M.L.~Gandarias for sending him her
related work.

\label{sophocleous-lp}

\end{document}